\documentclass[a4paper,oneside,10pt]{article}

\usepackage[english]{babel}    % slovenian language and hyphenation
\usepackage[utf8]{inputenc}    % make čšž work on input
\usepackage[T1]{fontenc}       % make čšž work on output
\usepackage{amsmath}    % basic ams math environments and symbols
\usepackage{amssymb,amsthm}    % ams symbols and theorems
\usepackage{mathtools}         % extends ams with arrows and stuff
\usepackage{url}               % \url and \href for links
\usepackage{icomma}            % make comma a thousands separator with correct 
\usepackage{units}             % \unit[1]{m} and unitfrac
\usepackage{enumerate}         % enumerate style
\usepackage{array}             % mutirow
\usepackage[usenames]{color}   % colors with names
\usepackage{graphicx}          % images
\usepackage{caption}           % captions
\usepackage{times}             % font
\usepackage{titlesec}          % section titles
\usepackage{soul}
\usepackage[ruled]{algorithm2e}% algorithm
\usepackage{todonotes}         % todo
\usepackage[bookmarks, colorlinks=true, linkcolor=black, anchorcolor=black,
citecolor=black, filecolor=black, menucolor=black, runcolor=black,
urlcolor=black, pdfencoding=unicode]{hyperref}  % clickable references, pdf toc
\usepackage[
  paper=a4paper,
  top=2.5cm,
  bottom=2.5cm,
left=2cm,
right=2cm
]{geometry}  % page geomerty

\newcommand{\R}{\mathbb{R}}
\newcommand{\T}{\mathsf{T}}
\renewcommand{\L}{\mathcal{L}}
\renewcommand{\b}{\boldsymbol}

\newcommand{\x}{\b{x}}
\newcommand{\w}{\b{w}}

\newcommand{\lap}{\nabla^2}

\newenvironment{itemize*}{\vspace{-6pt}\begin{itemize}\setlength{\itemsep}{0pt}\setlength{\parskip}{2pt}}{\end{itemize}}
\newenvironment{enumerate*}{\vspace{-6pt}\begin{enumerate}\setlength{\itemsep}{0pt}\setlength{\parskip}{2pt}}{\end{enumerate}}
\newenvironment{description*}{\vspace{-6pt}\begin{description}\setlength{\itemsep}{0pt}\setlength{\parskip}{2pt}}{\end{description}}

\newcommand{\Title}{Stability analysis of RBF-FD and WLS based local strong form meshless methods on scattered nodes}
\newcommand{\Author}{Mitja Jančič}
\newcommand{\AuthorTwo}{Gregor Kosec}

\title{\Title}
\author{\Author}
\date{\today}
\hypersetup{pdftitle={\Title}, pdfauthor={\Author}, pdfcreator={\Author},
pdfproducer={\Author}, pdfsubject={}, pdfkeywords={}}  % setup pdf metadata

\usepackage{multicol}
\usepackage{nth}
\usepackage[noabbrev]{cleveref}

\newenvironment{Figure}
{\par\medskip\noindent\minipage{\linewidth}}
{\endminipage\par\medskip}

\titleformat{\section}{\scshape\normalsize\centering}{\Roman{section}.}{1em}{}
\titleformat{\subsection}{\itshape\normalsize\centering}{\Alph{subsection}.}{1em}{}
\titleformat{\subsubsection}{\itshape\normalsize}{\arabic{subsubsection})}{1em}{}

\titlespacing\section{0pt}{2.82mm plus 2mm minus 0mm}{-0.71mm}
\titlespacing\subsection{5.5mm}{2.12mm}{-1.06mm}
\titlespacing\subsubsection{3.18mm}{0.72mm}{-2.12mm}
\titlespacing\paragraph{5.08mm}{0.0mm}{1ex}
\begin{document}

% NASLOV and stuff, ne spreminjaj.

\begin{center}
  \fontsize{24}{28}\selectfont
  \Title \\[3ex]
  \fontsize{11}{11}\selectfont
  \Author$^{1, 2}$,  \AuthorTwo$^{1}$ \\[0.71mm]
  \fontsize{10}{10}\selectfont
  $^1$ ``Jožef Stefan'' Institute, Parallel and Distributed Systems
  Laboratory, Ljubljana, Slovenia \\[0.5mm]
  $^2$ ``Jožef Stefan'' International Postgraduate School, Ljubljana,
  Slovenia \\[1mm]
  
  %%%%%%%%%%%%%%%%%%%%%%%%%%%%%%%
  %%%% NASTAVI PODATKE mail %%%%%
  %%%%%%%%%%%%%%%%%%%%%%%%%%%%%%%
  \href{mailto:Mitja.Jancic@ijs.si}{mitja.jancic@ijs.si},
  \href{mailto:Gregor.Kosec@ijs.si}{gregor.kosec@ijs.si}
\end{center}

\vspace{1ex}

\begin{multicols}{2}
  
  \fontsize{9}{9}\selectfont
  {\bfseries\noindent

  \emph{Abstract} -- The popularity of local meshless methods in the field of numerical simulations has increased greatly in recent years. This is mainly due to the fact that they can operate on scattered nodes and that they allow a direct control over the approximation order and basis functions. In this paper we analyse two popular variants of local strong form meshless methods, namely the radial basis function-generated finite differences (RBF-FD) using polyharmonic splines (PHS) augmented with monomials, and the weighted least squares (WLS) approach using only monomials. Our analysis focuses on the accuracy and stability of the numerical solution computed on scattered nodes in a two- and three-dimensional domain. We show that while the WLS variant is a better choice when lower order approximations are sufficient, the RBF-FD variant exhibits a more stable behavior and a higher accuracy of the numerical solution for higher order approximations, but at the cost of higher computational complexity.
  }
  \fontsize{10}{11}\selectfont
  
  \bigskip
  {\bfseries \fontsize{9}{9}
    \emph{Keywords -- meshless; WLS; RBF-FD; stability; scattered nodes}
  }
  
  %%%%%%%%%%%%%%%%%%%%%%%%%%%%%%%
  %%%% Chapters/sections %%%%%%%%
  %%%%%%%%%%%%%%%%%%%%%%%%%%%%%%%

  \section{Introduction}
  \label{sec:introduction}
  Computational science has become an important aspect of technological advancement in the fields of science and engineering. Thanks to the unprecedented computing power at our disposal, many real-life problems are being numerically treated to deepen our understanding of a phenomenon under consideration.

  In the field of numerical simulations, meshless methods are becoming increasingly popular with recent uses in the fields of fluid mechanics~\cite{rot2022refined}, linear elasticity~\cite{slak2019refined}, contact problems~\cite{slak2019adaptive}, advection-dominated problems~\cite{mavrivc2020equivalent} and even in financial sector~\cite{milovanovic2018radial}. Historically, mesh-free methods were introduced in the 1970s with the smoothed particle hydrodynamics (SPH)~\cite{lucy1977numerical, gingold1977smoothed} and then followed by several generalizations of the Finite Difference Method (FDM), e.g.\ the Finite Point Method~\cite{onate2001finite}, the Generalized Finite Difference Method~\cite{gavete2003improvements} and the Radial Basis Function-Generated Finite Differences (RBF-FD)~\cite{tolstykh2003using}. Nowadays, a lot of research is also devoted to reducing the computational time by employing the advantages of modern computer architecture~\cite{DepolliTrobec2019, TrobecDepolli2021}.

  The ability of meshless methods to operate on scattered nodes makes them very attractive in many real-life cases, where the domain shapes are often non-trivial. This is mainly because node positioning is easier than mesh generation (required by the mesh-based methods). Several algorithms for node positioning have been proposed to the meshless community. Some even support variable node density distributions~\cite{slak2019generation} and employing parallelization~\cite{depolli2022parallel} to reduce computational time. Another attractive feature of meshless methods is that the linear differential operator approximation allows a direct control over the order of the approximation method, as demonstrated in~\cite{janvcivc2021monomial, janvcivc2021p}, which can effectively be used to increase stability, as will be shown in this work.

  With many proposed meshless variants it is often not clear which is the most optimal in terms of numerical stability. Therefore, the aim of this paper is to compare the two commonly used variants, namely the WLS with monomials, also known as diffuse approximation method, and the RBF-FD variant with polyharmonic splines (PHS) augmented with monomials. The stability of the two methods is evaluated by solving a Poisson problem in two- and three-dimensional domain for lower and higher order approximations.

  The rest of the paper is organized as follows: In Section~\ref{sec:meshless} both WLS and RBF-FD approximation methods are presented, in Section~\ref{sec:problem} our case study is presented, in Section~\ref{sec:results} the results are shown and commented. Finally, in Section~\ref{sec:conclusions} conclusions and our findings are given.

  \section{Meshless methods}
  \label{sec:meshless}
  In meshless methods, a linear differential operator $\L$ at each node $\b x_c$ from the domain space $\Omega$ is approximated over a set of nearby nodes
  \begin{equation}
      \label{eq:ansatz}
      \widehat{\L u}(\b x_c) \approx \sum_{i = 1}^n w_iu(\b x_i) =\b w_\L \b u
  \end{equation}
  for any function $u$, $n$ nearby nodes also known as \emph{stencil nodes} and weights $w_i$, that are obtained by enforcing the equality of the equation~\eqref{eq:ansatz} for a given set of $s$ basis functions $\left\{ p _j \right\} _{j=1}^s$. The most common choice is to declare the nearest $n$ nodes as stencil, but some authors reported special stencil selection algorithms that increase the overall stability of the approximation~\cite{jacquemin2021unified, davydov2011adaptive}.

  The approximation~\eqref{eq:ansatz} is general in the sense that it holds for any linear differential operator $\L$ for any support size $n$ and any type or number $s$ of chosen basis functions $p$. As long as the number of basis functions is equal to the number of support nodes ($s=n$), the formulation of~\eqref{eq:ansatz} yields a quadratic system of equations
  \begin{equation}
      \label{eq:system}
      \underbrace{\begin{bmatrix}
                      p _1(\b x_1) & \cdots & p _1 (\b x_n) \\
                      \vdots       & \ddots & \vdots        \\
                      p _s(\b x_1) & \cdots & p_s(\b x_n)
      \end{bmatrix}}_{\b{\mathrm{P}}}
      \underbrace{\begin{bmatrix}
                      w_1    \\
                      \vdots \\
                      w_n
      \end{bmatrix}}_{\b w}=
      \underbrace{\begin{bmatrix}
      (\L p _1)(\b x_c)
                      \\
                      \vdots                  \\
                      (\L p _s)(\b x_c)
      \end{bmatrix}}_{\b l}
  \end{equation}
  as is the case with the so-called local collocation methods~\cite{kosec2008local}. However, larger support sizes are often used, resulting in an overdetermined system of equations. In such cases, the linear system is usually treated as a minimization of the weighted least squares (WLS) norm~\cite{sirca2012computational}
  \begin{align}
      \label{eq:wls-error}
      \begin{split}
          \left\| e \right\|_{2, w}&=\sqrt{\sum_{i=1}^n w^2 ( \widehat{u} (\b x_i) - u_i)^2} =\\ &=\sqrt{\sum_{i=1}^2(we_i)^2}.
      \end{split}
  \end{align}

  Common basis functions include: Multiquadrics, Gaussians, Radial Basis Functions (RBFs) and Monomials. In this paper we focus on two different types of basis functions, i.e., monomials and polyharmonic splines (PHS) augmented with monomials, resulting in two variants of meshless methods also known as the WLS approach and the RBF-FD variant respectively. While WLS approximation using a set $\left\{ p _j \right\} _{j=1}^s$ monomials up to and including degree $m$ as basis functions is fully defined above, the RBF-FD approximation is defined in the following section.

  \subsection{Radial Basis Function-Generated Finite Differences}
  \label{sec:rbffd}
  Let us take RBFs $\varphi$, such that $\varphi : [0,\infty) \rightarrow \R$ is centered at the stencil nodes of a central node $\b x_c$. The matrix $\b \Phi$ from the linear system~\eqref{eq:system} is then obtained by evaluating basis functions
  \begin{equation}
      \Phi _{ij} = \varphi(\left\| \b x_i - \b x_j \right\|)
  \end{equation}
  and the vector $\b l$ is assembled by applying the considered operator $\L$ to the basis functions evaluated at $\b x_c$, i.e.,
  \begin{equation}
      l_\varphi^i= (\L \varphi (\left\| \b x - \b x_i \right\|))\big|_{ \b x = \b x_c}.
  \end{equation}

  We can choose from different types of RBFs. Until recently, Hardy's multiquadrics or Gaussians were commonly used, but both depend on a shape parameter which effectively governs the accuracy and stability of the approximation~\cite{bayona2017role}. To avoid the shape parameter, we use PHS, defined as
  \begin{equation}
      \varphi(r) = \begin{cases}
                       r^k,       & k \text{ odd}  \\
                       r^k\log r, & k \text{ even}
      \end{cases},
  \end{equation}
  where $r$ denotes the Eucledian distance between two nodes.

  However, the use of pure RBFs as basis functions may lead to stagnation errors~\cite{flyer2016role}. Therefore, in addition to the RBFs, augmentation with monomials up to and including degree $m$ is added. This essentially means that we take a set of polynomials $\left\{ p _j \right\} _{j=1}^s$  with up to and including degree $m$ with $s=\binom{m+d}{d}$ and in addition to the RBF part of the approximation, the following exactness constraint for monomials is enforced
  \begin{equation}
      \label{eq:p-constraint}
      \sum_{i=1}^{s} w_i p_j (\x_i) = (\L p_j)(\x_c).
  \end{equation}
  The additional constraints make the approximation overdetermined,
  which is treated as a constrained optimization problem~\cite{flyer2016role}. For practical computation, the optimal solution can be expressed as a solution of
  a linear system
  \begin{equation}
      \label{eq:rbf-system-aug}
      \begin{bmatrix}
          \b{\b \Phi}       & \b{\mathrm{P}} \\
          \b{\mathrm{P}}^\T & 0
      \end{bmatrix}
      \begin{bmatrix}
          \w \\ \b \lambda
      \end{bmatrix}
      =
      \begin{bmatrix}
          \b \ell_{\varphi} \\ \b \ell_{p}
      \end{bmatrix},
  \end{equation}
  where $\b{\mathrm{P}}$ is an $n \times s$ matrix of polynomials evaluated at stencil nodes as is defined in a purely monomial approximation~\eqref{eq:system}, and $\b \ell_p$ is the vector of values assembled by applying the considered operator $\L$ to
  the polynomials at $\x_c$
  \begin{equation}
      l_p^i = (\L p_i(\x))\big|_{ \b x = \b x_c}
  \end{equation}
  and $\b \lambda$ are Lagrange multipliers.

  Finally, the system~\eqref{eq:rbf-system-aug} is solved to obtain weights and the approximate operator $\L$ at $\x_c$. Lagrange multipliers are discarded.

  Note that the exactness of~\eqref{eq:p-constraint} ensures the convergence behavior and also provides direct control over the convergence rate of the RBF-FD variant, since the local approximation has the same order as the polynomial basis used~\cite{bayona2017role}. Also notice that the linear system~\eqref{eq:rbf-system-aug} is larger when RBFs are augmented with monomials compared to when only monomials are used, making the RBF-FD approximation method computationally more expensive than the WLS approach.
  
  \subsection*{Implementation note}
  All elements of the solution procedure using meshless methods used in this paper are implemented in C++ using an object-oriented approach and a template system to achieve dimensionality independence. The numerical library used in this work and developed in-house is the \emph{Medusa library}~\cite{slak2021medusa}.
  
  \section{Case setup}
  \label{sec:problem}
  Any conclusions drawn from the analysis should not be specific to any particular domain shape or problem setup. We therefore choose a simple Poisson problem with Dirichlet boundary conditions on a $d$-dimensional sphere. All observations we make on this simple example should be understood as fundamental properties of the approximation method employed and therefore apply to all more complex domain shapes and problem setups.
  
  In the form of a PDE system, Poisson problem can be written as
  \begin{align}
    \label{eq:plap}
    \lap u(\x) & = f(\x) \qquad                 &  & \text{in } \Omega,         \\
    \qquad
    \label{eq:bc}
    u(\x)      & = \prod_{i=1}^d \sin(\pi x _i) &  & \text{on } \partial \Omega
  \end{align}
  where the right hand side was chosen to be
  \begin{equation}
    \label{eq:rhs}
    f(\x) = -d \pi  ^2 \prod _{i = 1}^d \sin(\pi x _i)
  \end{equation}
  and, as previously noted, the domain $\Omega$ is a $d$-dimensional sphere
  \begin{equation}
    \Omega = \left\{ \b r \in \R^d, \left\| \b r \right\|  \leq 1 \right\}.
  \end{equation}
  Example solution to~\cref{eq:plap,eq:bc,eq:rhs} is shown in Figure~\ref{fig:example-solution}. 
  
  \begin{Figure}
    \centering
    \includegraphics[width=\linewidth]{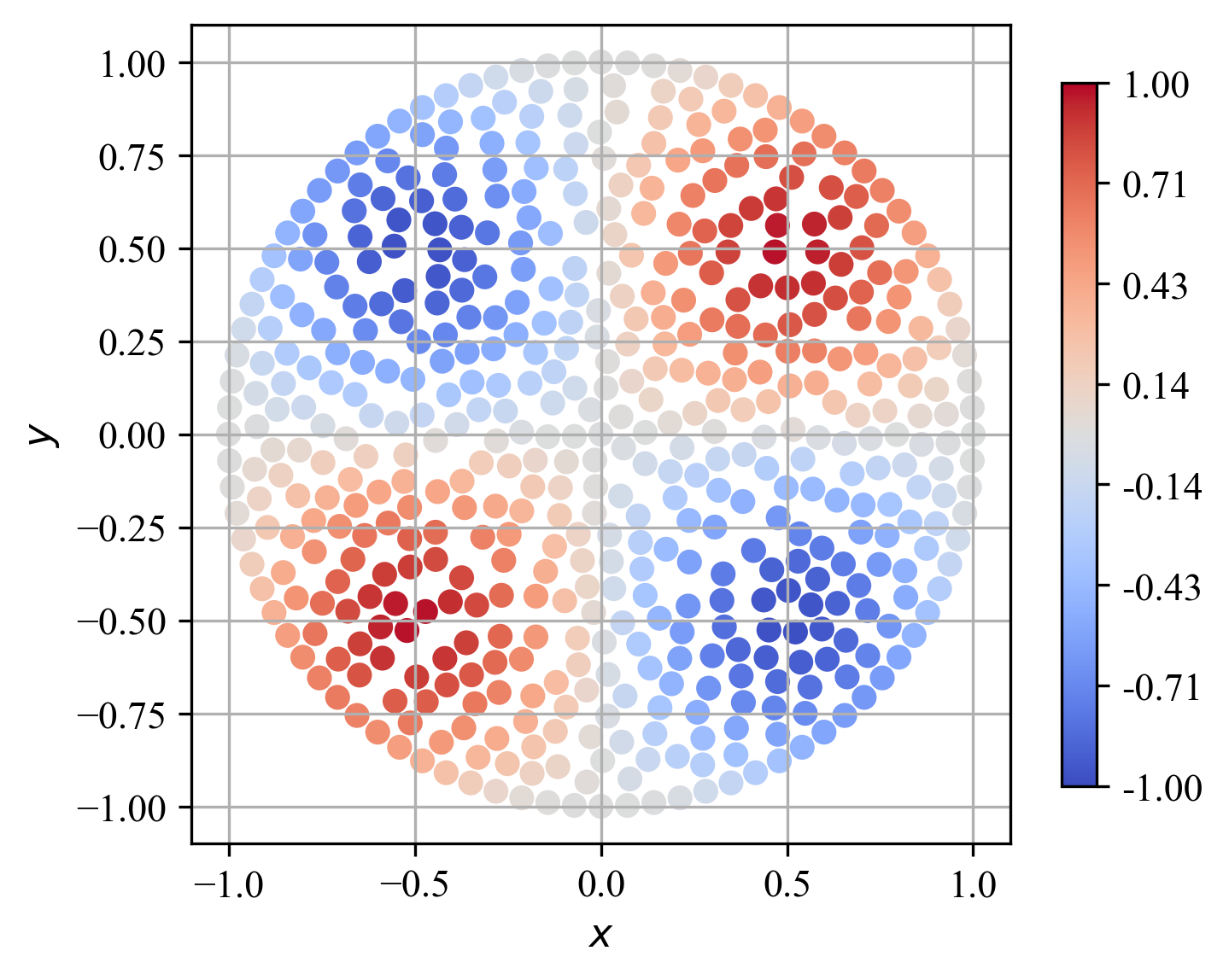}
    \captionof{figure}{Example solution to 2D Poisson problem with Dirichlet boundary conditions on $N=523$ scattered nodes.}
    \label{fig:example-solution}
  \end{Figure}
  
  The closed form solution of the above problem is $u(\x) = \prod_{i=1}^d \sin(\pi x _i)$. Having the closed form solution allows us to evaluate the accuracy of the numerical solution $\widehat{u}$ by computing the infinity norm error
  \begin{equation}
    e_\infty = \frac{\|\widehat{u} - u\|_\infty}{\|u\|_\infty}, \quad \|u\|_\infty = \max_{i=1, \ldots, N} |u_i|.
  \end{equation}
  The infinity norm was chosen because the authors have shown in~\cite{janvcivc2021monomial} that it measures the lowest convergence rates and does not involve averaging, unlike the commonly chosen $2$-norm error. 
  
  Numerical results are obtained using the two meshless methods described in Section~\ref{sec:meshless}, for a given domain discretization. Firstly, the RBF-FD using PHS radial basis function $\varphi(r) = r^3$ and monomial augmentation up to and including order $m \in \left\{ 2, 4, 6 \right\}$ is used, and secondly, the WLS approach using only monomials up to and including the same order $m$ is used. After both variants have been employed, the domain discretization is discarded and discretized again so that the number of nodes $N$ remains approximately the same, but the positions of the discretization nodes are different. The same process is repeated $N_{runs} = 100$ times. Note that after every domain discretization, the shapes must also be recomputed. This potentially leads to a different accuracy of the numerical solution and allows us to observe the stability of the numerical methods. 
  
  % For both approximation methods, the nearest $n$ nodes were selected as the stencil, where $n$ was defined as
  % \begin{equation}
  %   n = 2\binom{m+d}{d},
  % \end{equation}
  % as recommended by Bayona~\cite{bayona2017role} for the computationally more expensive RBF-FD approximations.
  
  The aim of this research is to determine which of the two approximation methods is more prone to a non-optimal discretization of the domain. For that, we introduce a normalized infinity norm error
  \begin{equation}
    \label{eq:spread}
    \frac{e_\infty ^{max} - e_\infty ^{min}}{e_\infty ^{median}}
  \end{equation}
  effectively describing the largest norm difference identified within the $N_{runs}=100$ runs, divided by the median value.
  
  \section{Results}
  \label{sec:results}
  In this section results are presented. All computations were performed with parallel execution on a computer with \texttt{AMD Threadripper 3990X} processor and \texttt{8x32GB DDR4} memory. The code\footnote{The source code is available at: \url{https://gitlab.com/e62Lab/public/cp-2022-mipro-engine_stability} under tag \emph{v1.2}.} was compiled using \texttt{g++ (GCC) 9.3.0 for Linux} with \texttt{-O3 -DNDEBUG -fopenmp} flags. 
  
  In all of the following figures, the blue colour is used to indicate the order of the approximation method $m=2$, red for $m=4$ and green for $m=6$. In the following subsections, the results for each dimension are presented separately, while our findings are summarized and presented as part of our conclusions in Section~\ref{sec:conclusions}. Only $d=\left\{ 2, 3 \right\}$ dimensional domains are studied, since scattered nodes do not make sense in $d=1$ dimensional problems and higher dimensional spaces $d > 3$ are beyond the scope of this paper.
  
  \subsection{The effect of stencil size}
  \label{sec:stencil_effect}
  This section presents the effect of stencil size on accuracy of the numerical solution and on the stability of the meshless variant. A scan over a range of stencil sizes $n$ is made and shown in Figure~\ref{fig:2d-scan}. Numerical solution has been obtained $N_{runs} = 100$ times at any given stencil size $n$, and with new domain discretization after every run keeping the total node count $N\approx 40600$. We can clearly observe that the error can be significantly higher if the support sizes are not sufficiently large, independent of the meshless variant. However, beyond that point, the dependency of the accuracy and the stability of the numerical solution on the stencil size is practically negligible. The tipping point for both approximation methods is in the neighborhood of recommendations made by Bayona~\cite{bayona2017role} for the RBF-FD, that is 
  \begin{equation}
    \label{eq:binom}
    n = 2\binom{m+d}{d}.
  \end{equation}

  \begin{Figure}
      \centering
      \includegraphics[width=\linewidth]{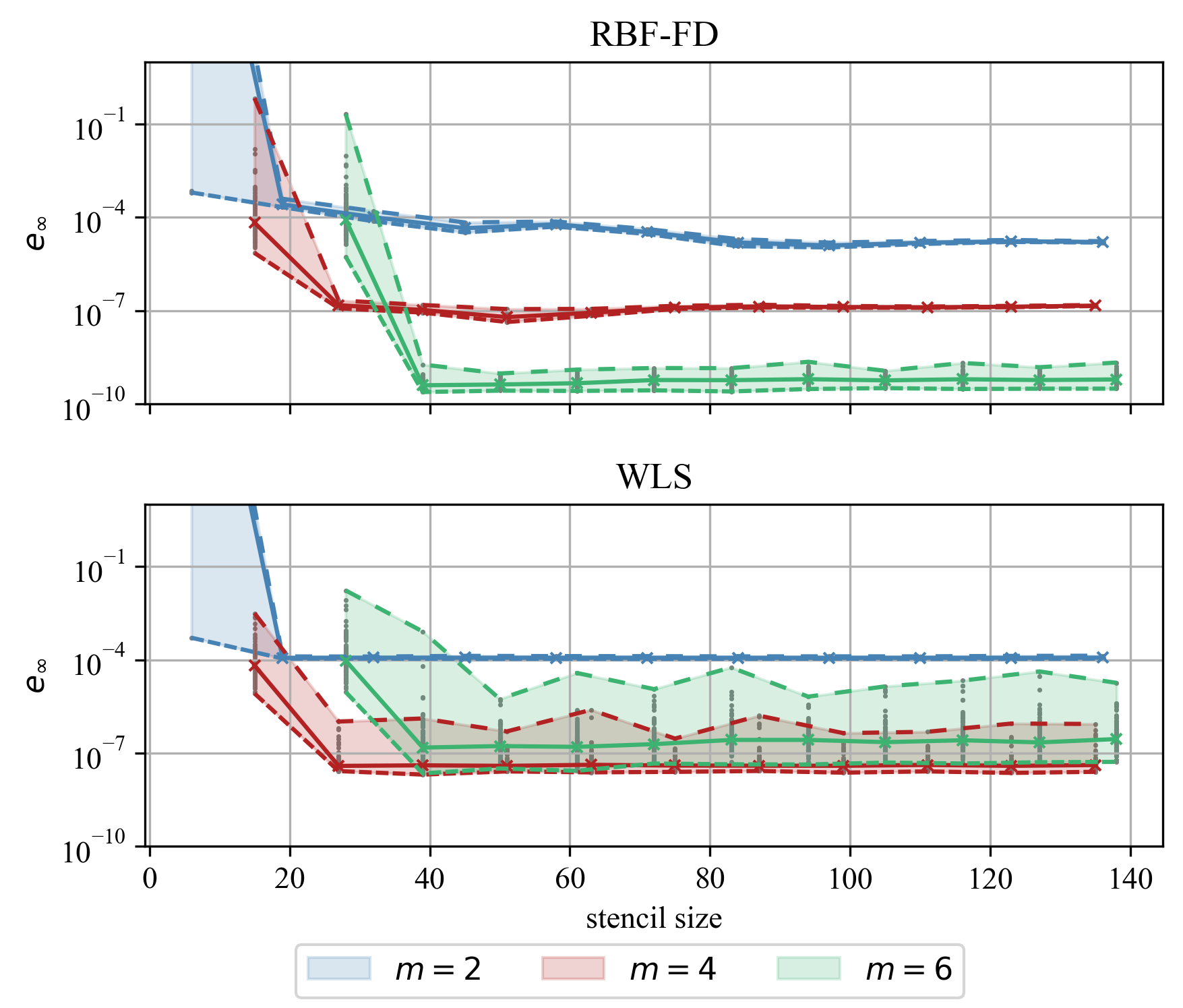}
      \captionof{figure}{Support size scan in two-dimensional domain with approximately $N\approx 40600$ nodes.}
      \label{fig:2d-scan}
  \end{Figure}

  Although all of the above conclusions were made on a 2-dimensional case, similar observations can be made in a three-dimensional domain, shown in Figure~\ref{fig:3d-scan}. However, in a three-dimensional case, a smaller number of discretization points $N\approx 28100$ has been used due to the fact that the support sizes are generally larger (see equation~\eqref{eq:binom}), making the computational times longer.

  \begin{Figure}
      \centering
      \includegraphics[width=\linewidth]{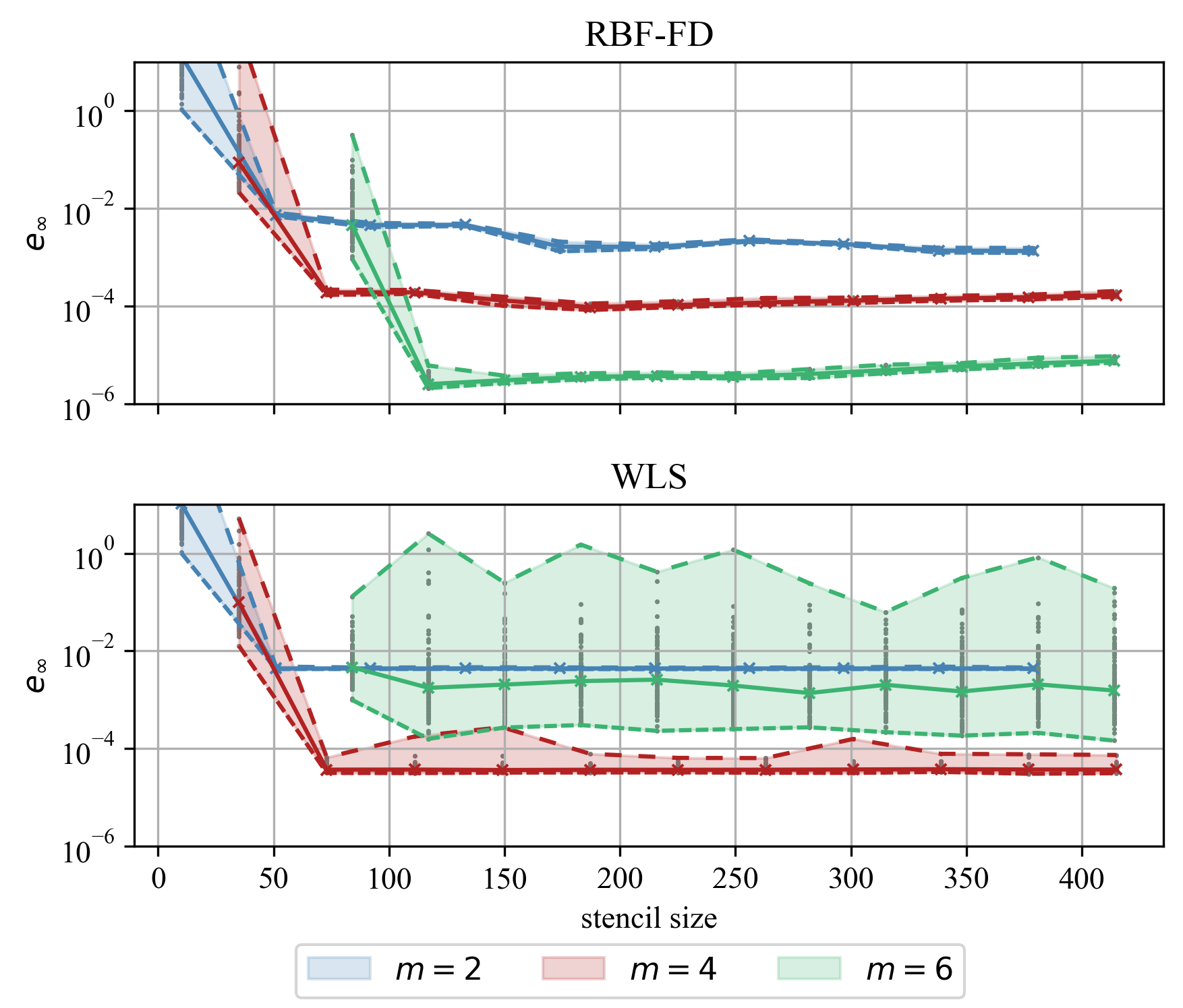}
      \captionof{figure}{Support size scan in three-dimensional domain with approximately $N\approx 28100$ nodes.}
      \label{fig:3d-scan}
  \end{Figure}

  \subsection{Convergence rates}
  \label{sec:convergence}
  % \subsubsection*{Two-dimensional domain}
  Considering the observations from the previous subsection, we continue our analysis by limiting ourselves to a single support size as defined in equation~\eqref{eq:binom}.

  Convergence rates in the case of a two-dimensional domain are shown in Figure~\ref{fig:2d-results}. The fact that the errors eventually diverge is a consequence of the errors in finite precision arithmetic, as previously observed by Flyer~\cite{flyer2016role}. As expected, the order of magnitude of the infinity norm error is the same for both approximation methods, but small differences can be observed. First, the accuracy achieved with the higher order ($m=6$) approximation is significantly better with the RBF-FD than with the WLS approach, and second, the spread around a median error value is significantly smaller for the RBF-FD.

  \begin{Figure}
      \centering
      \includegraphics[width=\linewidth]{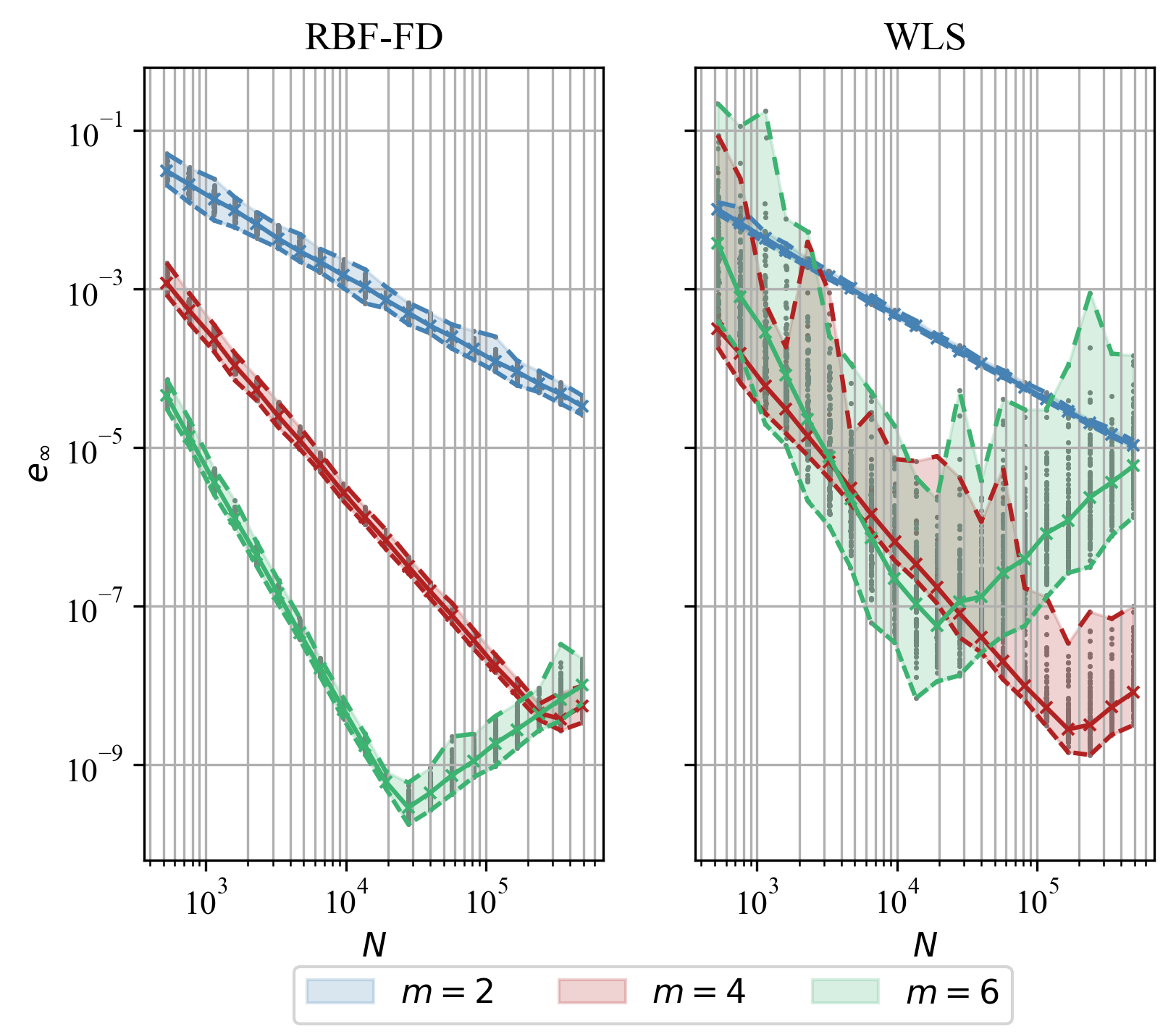}
      \captionof{figure}{Convergence rates in two-dimensional domain.}
      \label{fig:2d-results}
  \end{Figure}

  The spread observed after successive $N_{runs} = 100$ runs is further examined in Figure~\ref{fig:2d-spread}, where the normalized spread computed as defined in~\eqref{eq:spread} is evaluated and shown. We can see that the normalized spread is on average approximately constant, i.e., independent of the number of discretization points and approximately equal to 1 for a RBF-FD variant. We also find that the spread is about two orders of magnitude larger and unpredictable for the WLS approximation - with one exception, that is the low order WLS approximation ($m=2$), which clearly outperforms the RBF-FD variant in terms of stability and precision. In general, in two-dimensional domains, the RBF-FD variant is more stable, achieving better accuracy for higher order approximations, while lower order approximations are more stable and computationally cheaper to obtain with the WLS variant.
  
  \begin{Figure}
    \centering
    \includegraphics[width=\linewidth]{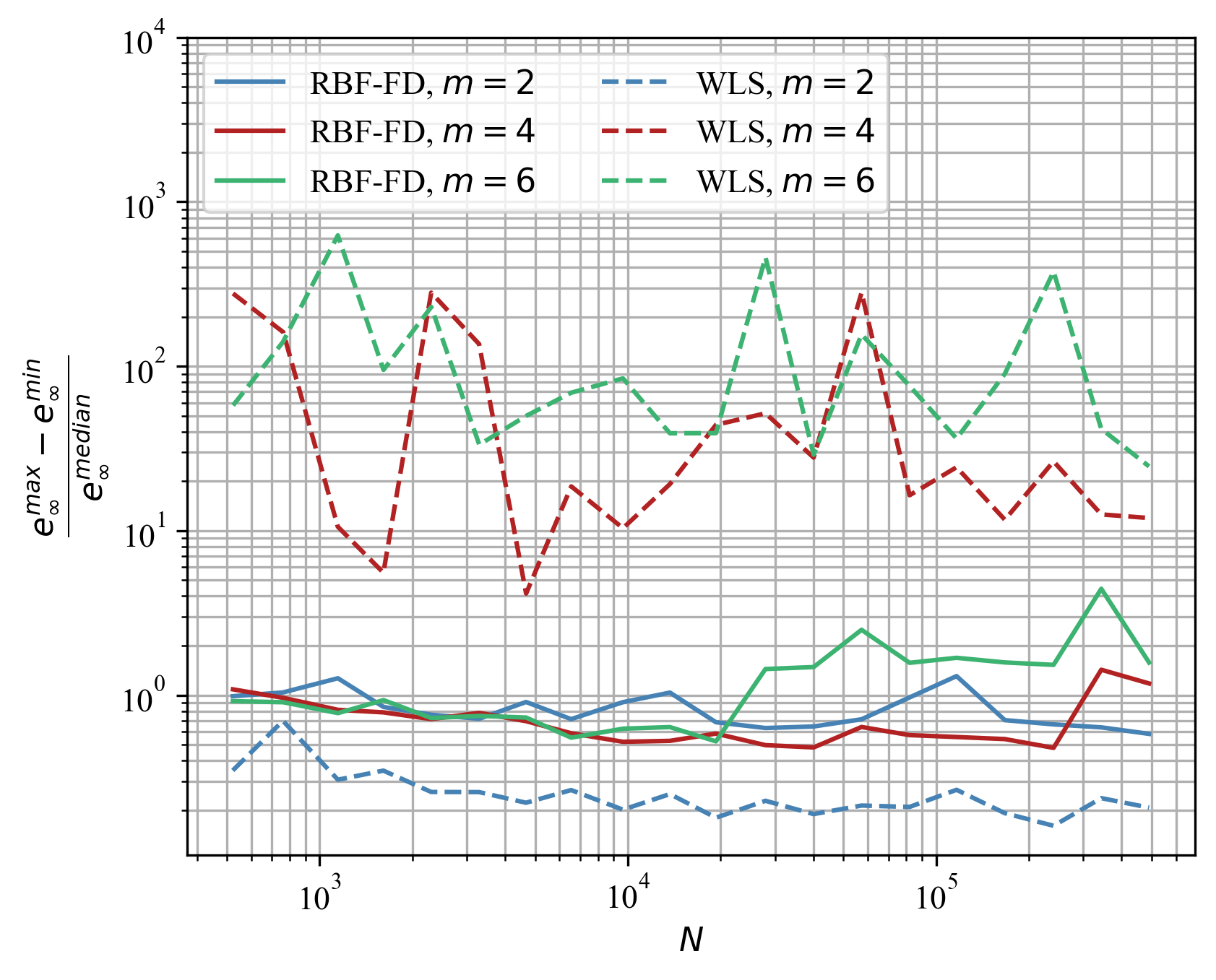}
    \captionof{figure}{Normalized spread size around median value in two-dimensional domain.}
    \label{fig:2d-spread}
  \end{Figure}
  
  Furthermore, convergence rates for a three-dimensional case are shown in Figure~\ref{fig:3d-results}. At first glance, we observe some similarities with the two-dimensional case in Figure~\ref{fig:2d-results}: Firstly, the errors are of the same order of magnitude for both meshless variants, and secondly, the spread size after successive $N_{runs} = 100$ runs is again in favour of the RBF-FD approximation method. More importantly, for the high order WLS approximation ($m=6$) and smaller number of discretization nodes $N\approx 10^3$, the infinity error norm is of the order of $10^1$. This essentially means that the WLS variant did not converge well. This is an important observation when studying stability, as we do not observe such a case in the results obtained with the RBF-FD variant.
  
  \begin{Figure}
    \centering
    \includegraphics[width=\linewidth]{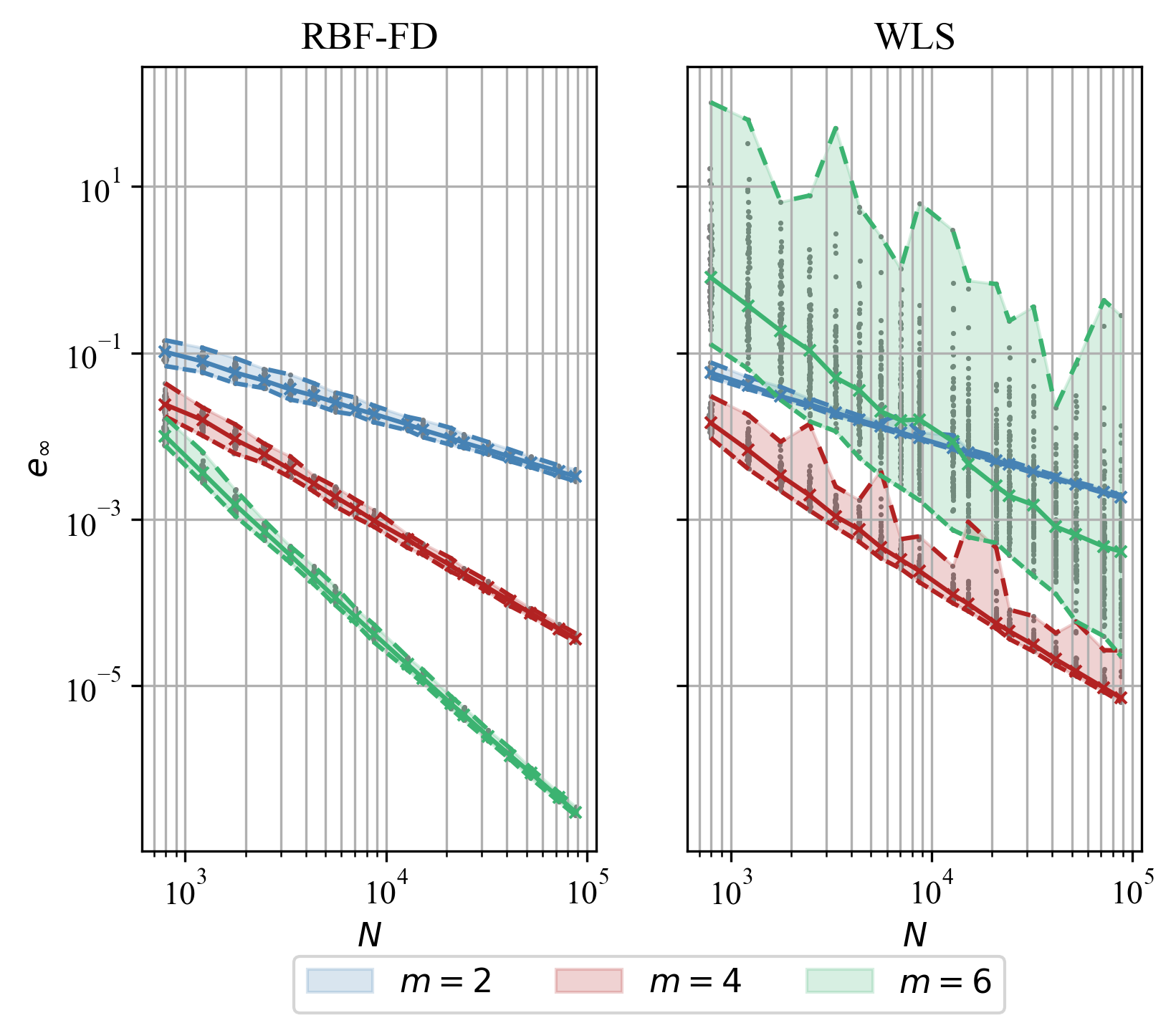}
    \captionof{figure}{Convergence rates in three-dimensional domain.}
    \label{fig:3d-results}
  \end{Figure}
  \begin{Figure}
    \centering
    \includegraphics[width=\linewidth]{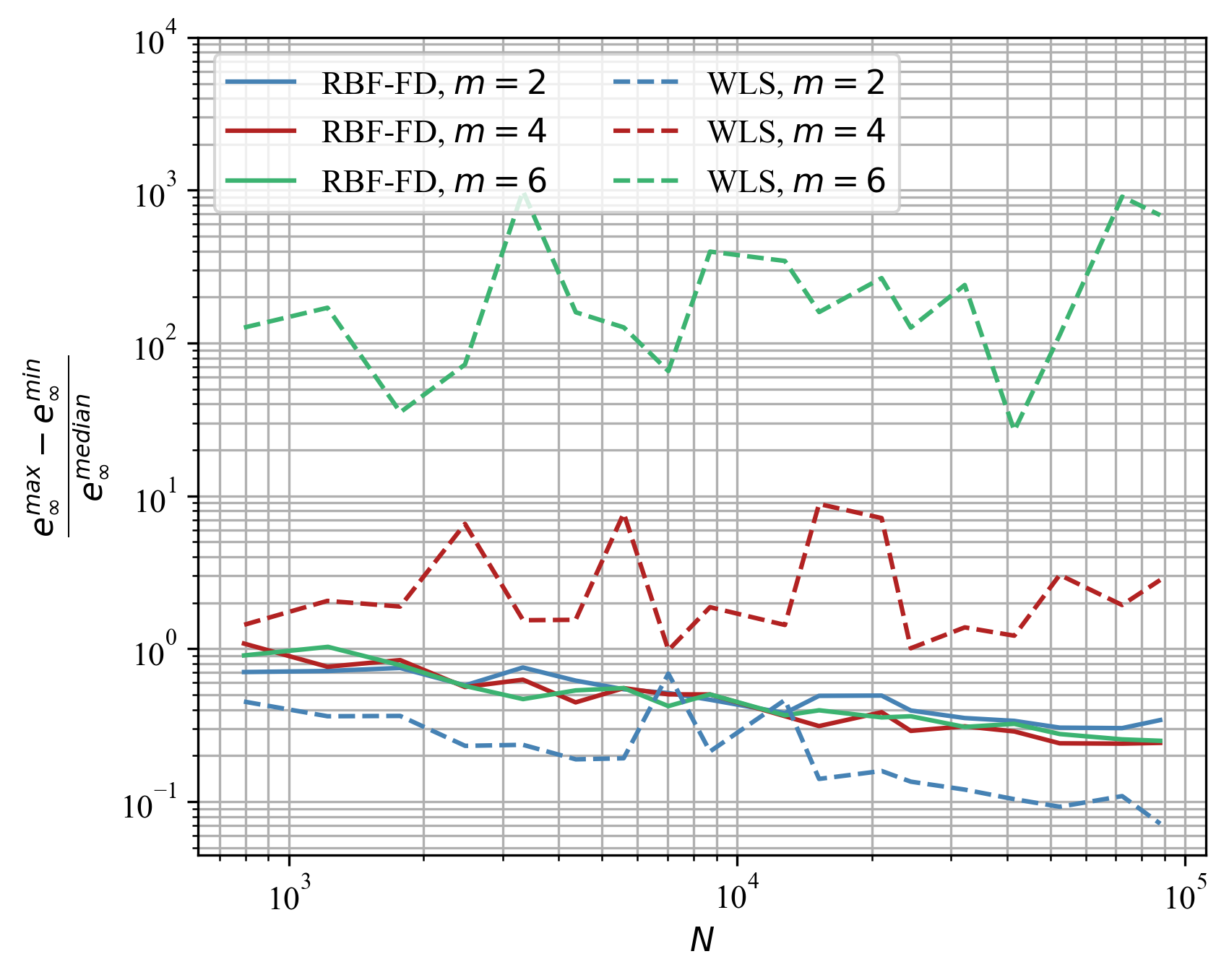}
    \captionof{figure}{Normalized spread size around median value in three-dimensional domain.}
    \label{fig:3d-spread}
  \end{Figure}

  Stability in a three-dimensional space is further studied in Figure~\ref{fig:3d-spread}. Similarly to the two-dimensional case in Figure~\ref{fig:2d-spread}, we find that the RBF-FD variant is substantially more stable than the WLS variant, especially for higher order approximations ($m > 2$). We also find that the amount of spread is slowly decreasing with the number of discretization points $N$ which was not observed in the two-dimensional case. This phenomenon is clearly seen for the RBF-FD approximations and for lower order WLS approximations ($m = 2$). The low order approximations ($m=2$) again appear to be in favour of the WLS meshless variant in terms of stability and computational complexity, although in some cases the normalized spread can be larger than that obtained by the RBF-FD (note the two peaks in the normalized spread for the low order WLS approximation in Figure~\ref{fig:3d-spread}). However, the higher order approximations are significantly more stable when obtained with the RBF-FD variant.

  \section{Conclusions}
  \label{sec:conclusions}
  In this paper, we compare the stability of two variants of meshless methods. We study solutions obtained with RBF-FD using polyharmonic splines augmented with monomials and with WLS approximation using only monomials as basis functions. Stability is assessed by solving a two- and three-dimensional Poisson problem with a tractable solution that allows us to evaluate the numerical solution in terms of the infinity norm error.

  We observe the effect of large enough stencil sizes has a negligible effect on the accuracy of the numerical solution and stability of the meshless variant. Additionally we show that in terms of stability, the RBF-FD variant can be several orders of magnitude better making it more prone to a non optimal domain discretization. This is particularly evident for the higher order approximations ($m > 2$), while the lower order approximations ($m=2$) are better and computationally cheaper to obtain using the WLS variant of meshless methods. We also find that the accuracy of the numerical solution obtained with a higher order approximation method can be significantly better when using RBF-FD than when using WLS.

  Further research is required to provide more accurate and descriptive guidelines as to which approximation method is most appropriate in particular cases. Although we provide some comments on the effect of stencil selection, we believe this aspect is in need of a more detailed study to make the methods not only stable but also computationally effective.

  According to our observations, RBF-FD variant is best used for problems that require a higher order approximation, while lower order approximations return better and faster results with WLS variant.

  \section*{Acknowledgments}
  \label{sec:ack}
  The authors would like to acknowledge the financial support of the ARRS research core funding No.\ P2-0095, ARRS project funding No.\ J2-3048 and the World Federation of Scientists.

  \bibliographystyle{unsrt}
  {\small
    \bibliography{ref}
  }
  
\end{multicols}

\end{document}